\documentclass[smallextended]{svjour3}       
\smartqed  
\usepackage{graphicx}
\usepackage[cmex10]{amsmath}
\usepackage{amssymb}
\usepackage{amsbsy}
\usepackage{algorithmicx,algorithm}
\usepackage{algpseudocode}
\usepackage{threeparttable}
\usepackage{enumerate}
 \usepackage{mathptmx}      
\newtheorem{conjec}{Conjecture}
\newtheorem{thm1}{Theorem}

\newtheorem{thm3}[thm1]{Theorem}
\newtheorem{thm4}[thm1]{Theorem}
\newtheorem{thm5}[thm1]{Theorem}
\newtheorem{lma1}{Lemma}
\newtheorem{lma2}[lma1]{Lemma}
\newtheorem{lma3}[lma1]{Lemma}

\begin{document}

\title{Extension of the Helleseth-Zinoviev result on the system of equations from Goethals codes and  Kloosterman sums
\thanks{Communicated by Pascale Charpin, Alexander Pott, Dieter Jungnickel}
}


\author{Minglong Qi \and Shengwu Xiong  }


\institute{Minglong Qi \and Shengwu Xiong  \at
              School of Computer Science and Technology,  Wuhan University of Technology \\
              Mafangshan West Campus, 430070 Wuhan City, China\\
              \email{mlqiecully@163.com (Minglong Qi) }\\
              \email{xiongsw@whut.edu.cn (Shengwu Xiong)}\\ 
}

\date{Received: date / Accepted: date}

\maketitle

\begin{abstract}
In the paper of Tor Helleseth and Victor Zinoviev (Designs, Codes and Cryptography, \textbf{17}, 269-288(1999)), the number of solutions of the system of equations from $ Z_{4} $-linear Goethals codes $ G_{4} $ was determined and stated in Theorem 4.  We found that Theorem 4 is wrong for $ m $ even. In this note, we complete Theorem 4, and give new divisibility modulo 12 of the Kloosterman sums deduced from Theorem 4, which is different from the result of the same authors for $ K(a^{4}+a^{3}) $ modulo 12.

\end{abstract}

\section{Introduction}
Let $ m $ be a positive integer,  $q=2^{m}, F:=\mathbb{F}_{q} $ be the finite field of $ q $ elements, $ F^{*}:=F\setminus \lbrace 0\rbrace $,  and $ F^{**}:=F\setminus \lbrace 0, 1\rbrace $. The well known Kloosterman sums \cite{carlitz,charpin} are defined by 
\begin{equation}\label{kloosterman_sum}
K(a)=\sum\limits_{x\in F^{*}}(-1)^{\mathrm{Tr}(ax+1/x)}
\end{equation}
where $ a\in F^{*} $, and $ \mathrm{Tr}(\cdot) $ is the trace function of $ F $ over $ \mathbb{F}_{2} $.

Let $ b,c\in F $. The problem of finding the coset weight distribution of $ Z_{4} $-linear Goethals codes $ G_{4} $,  is transformed into solving the following nonlinear system of equations over $ F $ \cite{tor-victor}:
\begin{equation}\label{G4_system}
\begin{cases}
x+y+z+u &=1\\
u^{2}+xy+xz+xu+yz+yu+zu &=b^{2}\\
x^{3}+y^{3}+z^{3}+u^{3} &=c
\end{cases}
\end{equation}
where $ x, y, z $ and $ u $ are pairwise distinct elements of $ F $. The number of solutions of (\ref{G4_system}) (see p.284 of \cite{tor-victor}), denoted by $ \mu_{2}(b,c) $, is given by 
\begin{equation}\label{mu_formula}
\mu_{2}(b,c)=\dfrac{2}{3}\times
\begin{cases}
M_{2}(b,c) &\qquad\text{if}\ \mathrm{Tr}(c)\ne \mathrm{Tr}(1),\\
M_{2}(b,c)-1 &\qquad\text{if}\ \mathrm{Tr}(c)= \mathrm{Tr}(1).
\end{cases}
\end{equation}
The $ M_{2}(b,c) $ in (\ref{mu_formula}) is defined by 
\begin{equation}\label{M2_formula}
4M_{2}(b,c)=\sum\limits_{v\in F^{**}}\bigl( 1+(-1)^{\mathrm{Tr}(\frac{k_{1}}{v})+l}\bigr) \bigl( 1+(-1)^{\mathrm{Tr}(\frac{k_{2}}{v})+l}\bigr)
\end{equation}
where $ k_{1}=b^{2}+c+1,  k_{2}=b^{2}+b+c+\sqrt{c} $,  and $ l=\mathrm{Tr}(b) $.
The explicit evaluation of $ M_{2}(b,c) $ is given by the following equation (see  (33) of \cite{tor-victor}):
\begin{equation}\label{M2_explicit}
M_{2}(b,c)=\dfrac{1}{4}\bigl(q-3+(-1)^{\mathrm{Tr}(k_{1})}K(k_{1}k_{2})-(-1)^{\mathrm{Tr}(b)}(3+(-1)^{\mathrm{Tr}(k_{1})})\bigr).
\end{equation}

We found that Theorem 4 of \cite{tor-victor}, which gives the explicit evaluation of $ \mu_{2}(b,c) $, is wrong for  $ m $ even. It is obvious that the authors of \cite{tor-victor} forgot to take account of the fact, that $ \mathrm{Tr}(1)=1 $ for $ m $ odd and 0 for $ m $ even, in the last step of the proof of Theorem 4. The following is the correct version of Theorem 4 of \cite{tor-victor}:
\begin{thm1}\label{thm1}
Let $ \mu_{2}(b,c) $ be the number of different 4-tuples $ (x, y, z, u) $, where $ x, y, z, u $ are pairwise distinct elements of $ F $, which are solutions to the system (\ref{G4_system}) over $ F $, where $ b, c $ are arbitrary elements of $ F $. 
\begin{enumerate}[(1)]
\item If $ m $ is odd and $ \mathrm{Tr}(c)=1 $ or $ m $ is even and $ \mathrm{Tr}(c)=0 $, then
\begin{equation*}
\mu_{2}(b,c)=\frac{1}{6}\bigl(q-8+(-1)^{\mathrm{Tr}(b)}(K(k_{1}k_{2})-3)\bigr).
\end{equation*}
\item If $ m $ is odd and $ \mathrm{Tr}(c)=0 $ or $ m $ is even and $ \mathrm{Tr}(c)=1 $, then
\begin{equation*}
\mu_{2}(b,c)=\frac{1}{6}\bigl(q-2-(-1)^{\mathrm{Tr}(b)}(K(k_{1}k_{2})+3)\bigr).
\end{equation*}
Where $  k_{1}=b^{2}+c+1$ and  $ k_{2}=b^{2}+b+c+\sqrt{c} $.
\end{enumerate}
\end{thm1}

From Theorem 4 of \cite{tor-victor}, the divisibility modulo 12 of the Kloosterman sums  $ K(a^{4}+a^{3}) $ is obtained and stated in Theorem 5 of \cite{tor-victor}. 
In Section 2, we will establish new divisibility modulo 12 of $ K(a^{4}+a^{3}) $.

\section{New divisibility modulo 12 of $ K(a^{4}+a^{3}) $}

Before stating the new results on $  K(a^{4}+a^{3}) $ modulo 12, we need several auxiliary results.
\begin{lma1}[\cite{carlitz}]\label{lemma1}
Let $ a\in F^{*} $. Then, $ K(a) =K(a^{2})$.
\end{lma1}

\begin{lma2}[\cite{lidl-nieder}]\label{lemma2}
The quadratic equation $ x^{2}+ax+b=0 $, where $ a\in F^{*} $ and $ b\in F $, has two distinct roots in $ F $, if and only if $ \mathrm{Tr}(b/a^{2})=0 $.
\end{lma2}

\begin{lma3}[\cite{lidl-nieder}]\label{lemma3}
Let $ a, b\in F $. Then,
\begin{enumerate}[(1)]
\item $ \mathrm{Tr}(a)=\mathrm{Tr}(a^{2}),  \mathrm{Tr}(a+b)=\mathrm{Tr}(a)+\mathrm{Tr}(b)$.
\item $ \mathrm{Tr}(1)=1 $ if $ m $ is odd, and $ \mathrm{Tr}(1)=0 $ if $ m $ is even.
\end{enumerate}
\end{lma3}

The following theorem shows that $  K(a^{4}+a^{3}) $ modulo 12 may have different values than Theorem 5 of \cite{tor-victor}  if $ a $ is under certain restriction:
\begin{thm3}\label{thm3}
Let $ a=b^{2}+b+1 $ where $ b\in F $. If $ a\ne 0 $, then 
\begin{equation*}
K(a^{4}+a^{3})\equiv
\begin{cases}
3 &\pmod {12},\ \text{if}\ m\equiv 1\pmod 2,\\
7 &\pmod {12},\ \text{if}\ m\equiv 0\pmod 2.
\end{cases}
\end{equation*}
\end{thm3}
\begin{proof}
  Set $ k_{2}=b^{2}+\sqrt{c} $, i.e., $ b^{2}+b+c+\sqrt{c}=b^{2}+\sqrt{c} $. We obtain $ b=c $, and $ k_{1}=b^{2}+b+1 $. Further, we have $ \mathrm{Tr}(b)=\mathrm{Tr}(c), \mathrm{Tr}(k_{1})=\mathrm{Tr}(b^{2}+b+1)=\mathrm{Tr}(1) $, and $ k_{1}^{2}k_{2}^{2}=b^{8}+b^{6}+b^{4}+b^{5}+b^{3}+b$. Let $ a=b^{2}+b+1 $. It is easy to see that $ k_{1}^{2}k_{2}^{2}=b^{8}+b^{6}+b^{4}+b^{5}+b^{3}+b=a^{4}+a^{3} $ in $ F $. By Lemma \ref{lemma1}, $ K(k_{1}k_{2})=K((k_{1}k_{2})^{2})=K(a^{4}+a^{3}) $. Substituting $ \mathrm{Tr}(k_{1})=\mathrm{Tr}(1) $ and $ \mathrm{Tr}(b)=\mathrm{Tr}(c) $ into (\ref{M2_explicit}),  using (\ref{mu_formula}), and taking account that $ \mathrm{Tr}(1)=1 $ if $ m $ is odd and 0 if $ m $ is even, we obtain that 
 \begin{enumerate}[(1)]
 \item Case that $ m $ is odd:
 \begin{equation*}
 \mu_{2}(b,c)=\frac{1}{6}\bigl(q-5-K(a^{4}+a^{3})\bigr)\ \text{for}\ \mathrm{Tr}(c)\in \lbrace 0, 1\rbrace.
 \end{equation*}
  \item Case that $ m $ is even:
  \begin{equation*}
  \mu_{2}(b,c)=
  \begin{cases}
  \frac{1}{6}\bigl(q-11+K(a^{4}+a^{3})\bigr), &\  \text{if}\ \mathrm{Tr}(c)=0,\\
  \frac{1}{6}\bigl(q+1+K(a^{4}+a^{3})\bigr), &\  \text{if}\ \mathrm{Tr}(c)=1.
  \end{cases}
  \end{equation*}
 \end{enumerate}
 
 The number $ 6\mu_{2}(b,c) $ is a multiple of 12 \cite{tor-victor}. It can be verified that $ (2^{m}-5)\equiv 3\pmod{12} $ if $ m $ is odd, and $ -(2^{m}-11)\equiv -(2^{m}+1) \equiv 7\pmod{12} $ if $ m $ is even. The theorem follows from  above argument.
\end{proof}\qed

Based on computer investigation, we propose the following conjecture:
\begin{conjec}\label{conjecture}
Let $ a=b^{2}+b+1 $ where $ b\in F $, and $ \xi=\sum_{i=0}^{n}a^{2i}+a$ with $ n\ge 0 $. If $ \xi\ne 0 $, then 
\begin{equation*}
K(\xi^{4}+\xi^{3})\equiv
\begin{cases}
3 &\pmod {12},\ \text{if}\ m\equiv 1\pmod 2,\\
7 &\pmod {12},\ \text{if}\ m\equiv 0\pmod 2.
\end{cases}
\end{equation*}
\end{conjec}

To illustrate Conjecture \ref{conjecture}, we present the following theorem:
\begin{thm4}
Let $ \xi=b^{4}+b+1 $ where $ b\in F $. If $ \xi\ne 0 $, then 
\begin{equation*}
K(\xi^{4}+\xi^{3})\equiv
\begin{cases}
3 &\pmod {12},\ \text{if}\ m\equiv 1\pmod 2,\\
7 &\pmod {12},\ \text{if}\ m\equiv 0\pmod 2.
\end{cases}
\end{equation*}
\end{thm4}
\begin{proof}
Let $ a=b^{2}+b+1 $ where $ b\in F $. Remark that $ \xi=a^{2}+a+1 $. Set $ k_{2}=b^{2}+b+c+\sqrt{c}=c+\sqrt{b} $. Then, $ c=b^{4}+b^{2}+b$ and $k_{1}=b^{4}+b+1 $. Further, we obtain that $ \mathrm{Tr}(c)=\mathrm{Tr}(b^{4}+b^{2}+b)=\mathrm{Tr}(b), \mathrm{Tr}(k_{1})=\mathrm{Tr}(b^{4}+b+1)=\mathrm{Tr}(1) $, and $ k_{1}^{2}k_{2}^{2}=(b^{8}+b^{2}+1)(b^{8}+b^{4}+b^{2}+b)=\xi^{4}+\xi^{3} $ in $ F $. The rest of the proof is analogous to Theorem \ref{thm3}, we omit the details.
\end{proof}\qed

If $ a\in F $ satisfies a certain condition, the binomial $ a^{4}+a^{3} $ of $ K(a^{4}+a^{3}) $ from Theorem 5 of \cite{tor-victor} may be reduced into a monomial of degree one, which is illustrated by the following theorem:
\begin{thm5}\label{thm5}
Let $ a\in F^{**} $ be a third root of unity in $ F^{*} $. Then,
\begin{equation*}
K(a)\equiv
\begin{cases}
7 &\pmod {12},\ \text{if}\ \mathrm{Tr}(a)=0,\\
11 &\pmod {12},\ \text{if}\ \mathrm{Tr}(a)=1.
\end{cases}
\end{equation*}
\end{thm5}
\begin{proof}
Set $ k_{1}=b^{2}+c+1=b+\sqrt{c} $, then $ k_{2}=b^{2}+b+c+\sqrt{c}=1 $. The condition $ k_{1}=b^{2}+c+1=b+\sqrt{c} $ holds if and only if $ \xi^{2}+\xi+1=0 $ where $ \xi= (b+\sqrt{c})\in F^{**}$, i.e., $ \xi $ is a third root of unity in $ F^{*} $. Further, we have $ K(k_{1}k_{2})=K(\xi) $. The rest of the proof is similar to Theorem \ref{thm3}, we omit the details.
\end{proof}\qed

\section{Conclusion}
In this note, we pointed out that Theorem 4 of \cite{tor-victor} is wrong for $ m $ even, and completed it by taking account of the fact that $ \mathrm{Tr}(1)=1 $ for $ m $ odd and 0 for $ m $ even in the proof of the theorem. We established some new divisibility modulo 12 of the Kloosterman sums $ K(a^{4}+a^{3}) $ too.


\begin{thebibliography}{99}
%
%
\bibitem{carlitz}
Carlitz, L.: Kloosterman sums and finite field extensions, Acta Arithmetica, \textbf{15}(2), 179-193 (1969).
\bibitem{charpin}
Charpin P., Helleseth T., Zinoviev V.: The divisibility modulo 24 of Kloosterman sums on $ GF(2^{m}), m $ odd, Journal of Combinatorial Theory, Series A, \textbf{144}, 322-338 (2007).
\bibitem{tor-victor}
Helleseth T., Zinoviev V.: On $ Z_{4} $-Linear Goethals Codes and Kloosterman Sums, Designs, Codes and Cryptography, \textbf{17}, 269-288 (1999).
\bibitem{lidl-nieder}
Lidl R., Niederreiter H.: Finite Fields, Encyclopedia of Mathematics and Its Applications, 2nd ed.,
Cambridge University Press, Vol. 20 (1997).
\end{thebibliography}


\end{document}